\documentclass[11pt,twoside,a4paper]{article} %\documentclass{article}
\usepackage{pgm} %available from http://mtr.utia.cas.cz/pgm06/PGMstyle/index.html
\usepackage{ifpdf}
\usepackage{amsthm}
\usepackage{amsmath}
\usepackage{amssymb}
\usepackage{graphicx}
\usepackage{psfrag}
\usepackage[all]{xypic}
\usepackage{url}

\newcommand{\sinpioverthree}{0.866025404}

\newcommand{\R}{\mathbb{R}}

\newcommand{\N}{\mathbb{N}}

\newcommand{\scrM}{\mathcal{M}}
\newcommand{\indep}{\! \perp \!\!\! \perp \!}

\newtheorem{theorem}{Theorem}
\newtheorem{lemma}[theorem]{Lemma}

\newtheorem{proposition}[theorem]{Proposition}
\newtheorem{observation}[theorem]{Observation}

\theoremstyle{definition}
\newtheorem{example}[theorem]{Example}%[section]
\newtheorem{remark}[theorem]{Remark}%[section]
\title{Geometry of Rank Tests}
\author{Jason Morton, Lior Pachter, Anne Shiu, Bernd Sturmfels, and
  Oliver Wienand
  \\Department of Mathematics, UC Berkeley}
\summary{We study partitions of the symmetric group
which have desirable geometric properties.
The statistical tests defined by such partitions involve counting all
permutations in the equivalence classes.
These permutations are the linear extensions of partially
ordered sets specified by the data.
Our methods refine rank tests of non-parametric statistics, such as
the sign test and the runs test, and are useful
for the exploratory analysis of ordinal data.
{\em Convex rank tests} correspond to probabilistic
conditional independence structures known
as semi-graphoids.
{\em Submodular rank tests} are classified by the faces of the cone
of submodular functions, or by Minkowski summands
of the permutohedron.
We enumerate all small instances of such rank tests.  
{\em Graphical tests} correspond to both graphical models
and to graph associahedra, and they have
 excellent statistical and algorithmic properties.}
\begin{document}
\maketitle

\section{Introduction}

The non-parametric approach to statistics was introduced by
\cite{Pitman1937SignificanceI}.
The emergence of microarray data in
molecular biology has led to a number of new tests for identifying significant
patterns in gene expression time series; see e.g.~\cite{Willbrand2005}.
This application motivated us to
develop a mathematical theory of rank tests.
We propose
that a {\em rank test} is a partition of $S_n$ induced by a
map $\, \tau : S_n \rightarrow T\,$
from the symmetric group of all permutations of $[n]=\{1,\ldots,n\}$
onto a set $T$ of statistics.
The statistic $\tau(\pi)$ is the {\em signature} of the permutation $\pi \in S_n$.
Each rank test defines a partition of $S_n$ into
classes, where $\pi$ and $\pi'$ are in the same class if and only if
$\tau (\pi) = \tau(\pi')$. We identify $T = {\rm image}(\tau)$ with
the set of all classes in this partition of $S_n$.
 Assuming the uniform distribution on  $ S_n$, the probability
of seeing a particular signature $t \in T$ is
$\,1/n! \,$ times $| \tau^{-1}( t)|$.
The computation of a $p$-value for a given permutation $\pi \in S_n$
typically amounts to summing
\begin{equation} \label{Pvalue}
{\rm Pr}(\pi') \quad = \quad
 \frac{1}{n !} \cdot |\, \tau^{-1} \bigl( \tau(\pi') \bigr)\, |
 \end{equation}
over all permutations $\pi'$ with
${\rm Pr}(\pi') < {\rm Pr}(\pi)$.
 In Section 2 we explain how existing rank tests
 can be understood from our  point of view.

In Section 3 we describe the class of {\em  convex rank tests} which captures
properties of tests used in practice. We work in the language
of algebraic combinatorics \cite{Stanley1997}.
 Convex rank tests are
 in bijection with polyhedral fans
that coarsen the hyperplane arrangement of $S_n$, and with
conditional independence structures known as
semi-graphoids \cite{Studeny2005Probabilistic}.

Section 4 is devoted to convex rank tests that are
induced by submodular functions.
These {\em submodular rank tests} are in
bijection with Minkowski summands
of the $(n-1)$-dimensional permutohedron and with structural
imset models. Furthermore, these tests are
at a suitable level of generality for the biological applications
that motivated us.
We make the connections to polytopes and independence models concrete
by classifying
all convex rank tests for $n \leq 5$.

In Section 5 we discuss the class of {\em graphical tests}.
In mathematics, these correspond to
graph associahedra, and in statistics
to graphical models.
The equivalence of these two structures is shown in
Theorem \ref{maingraphical}.
The implementation of convex rank tests requires the efficient enumeration of
linear extensions of partially ordered sets (posets).  A key ingredient is a
highly optimized method for computing distributive lattices.
Our software is discussed in Section~6.

\section{Rank tests and posets}

A permutation $\pi$ in $S_n$ is a
total order on $[n] = \{1,\ldots,n\}$.
This means that $\pi$ is a set
of $\binom{n}{2}$ ordered pairs
of elements in $[n]$.
If $\pi$ and $\pi'$ are permutations then
$\,\pi \cap \pi'\,$ is a partial order.

In the applications we have in mind, the data
are vectors $u \in \R^n$ with distinct coordinates.
The permutation associated with $u$ is the
total order $\,\pi = \{ \,(i,j)\in [n] \times [n] \,: \, u_i < u_j\,\}$.
We shall employ two other ways of writing this permutation.
The first is the {\em rank vector} $\,\rho = (\rho_1,\ldots,\rho_n)$,
whose defining properties are
$\{\rho_1,\ldots,\rho_n\} = [n]$ and
$\rho_i < \rho_j$ if and only if $u_i < u_j$.  That is, the coordinate of the rank vector with value $i$ is at the same position as the $i$th smallest coordinate of $u$.  
The second is the {\em descent vector}
$\delta = (\delta_1,\ldots,\delta_n)$, defined by $u_{\delta_i} > u_{\delta_{i+1}}$.  The $i$th coordinate of the descent vector is the position of the $i$th largest value of $u$.
For example, if $u = (11,7,13)$ then its permutation
is represented by
$\, \pi = \{ (2,1),(1,3),(2,3)\}$, by $\,\rho = (2,1,3)$, or by  $\,\delta = (3,1,2)$.

A permutation $\pi $ is a {\em linear extension} of a
partial order $P$ on $[n]$ if $P \subseteq \pi$.
 We write $\mathcal{L}(P) \subseteq S_n$ for the set of
linear extensions of~$P$.
A partition $\tau$ of the symmetric group $S_n$ is
a {\em pre-convex rank test}
if the following axiom holds:
\[\begin{array}{cc}
(PC) & \begin{array}{c} \text{If }\tau(\pi) = \tau(\pi') \text{ and } \pi'' \in \mathcal{L} (\pi \cap \pi') \\ \text{ then } \tau(\pi) \! = \! \tau(\pi') \!= \!\tau(\pi''). \end{array}\end{array}
\]
Note that $\, \pi'' \in \mathcal{L} (\pi \cap \pi') \,\,$  means
$\pi \cap \pi' \subseteq \pi''$.
For $n = 3$ the number of all rank tests is the Bell number $B_6 = 203$.
Of these $203$ set partitions of $S_3$, only $40$ satisfy the
axiom (PC).

Each class $C$ of a pre-convex rank test $\tau$ corresponds to
a poset $P$ on $[n]$; namely,  $P$ is the
intersection of all total orders in that class: $P=\bigcap_{\pi \in C} \pi$. The axiom
(PC) ensures that $C$ coincides with
the set $\mathcal{L}(P)$ of all linear extensions of $P$.  The inclusion $C \subseteq \mathcal{L}(P)$ is clear.  For the reverse inclusion, note that from any permutation $\pi$ in $ \mathcal{L}(P)$, we can obtain any other $\pi'$ in $ \mathcal{L}(P)$ by a sequence of reversals $(a,b) \mapsto (b,a)$,
where each intermediate $\hat{\pi}$ is also in $ \mathcal{L}(P)$.  Assume $\pi_0 \in \mathcal{L}(P)$ and $\pi_1 \in  C$ differ by one reversal  $(a,b)\in \pi_0$, $(b,a) \in \pi_1$.  Then
$(b,a) \notin P$, so there is some $\pi_2 \in C$
such that $(a,b) \in \pi_2$; thus, $\pi_0\in  \mathcal{L} (\pi_1 \cap \pi_2)$ by (PC).
This shows $\pi_0 \in C$.

A pre-convex rank test is thus
an unordered collection of posets $P_1,,\ldots,P_k$ on $[n]$ that
satisfies the property that  $S_n$ is the disjoint union of the subsets
$\mathcal{L}(P_1),\ldots, \mathcal{L}(P_k)$.
The posets $P_i$ that represent the classes in
a pre-convex rank test capture the
shapes of data vectors.

\begin{example}[The sign test for paired data] \label{sign_test}
The \emph{sign test} is performed on data that are paired as
two vectors $u=(u_1,u_2, \dots,u_m)$ and
$ v = ( v_1, v_2, \dots,  v_m)$.  The null hypothesis
is that the median of the differences $u_i - v_i$ is 0.
The test statistic is the number of differences
that are positive. This test is a rank test, because
$u$ and $v$ can be transformed into the
overall ranks of the $n=2m$ values, and the rank vector
entries can then be compared. This test coarsens the convex rank test which is
the MSS of Section 4 with $\mathcal{K}=\{\{1,m+1\},\{2, m+2\}, \dots \}$.
\end{example}

\begin{example}[Runs tests]
A \emph{runs test} can be used when there is a natural ordering on the data
points, such as in a time series. The data are
transformed into a sequence of `pluses' and `minuses,' and
the null hypothesis is that the number
of observed runs is
no more than that expected by chance.  A runs test is a coarsening of the convex rank test $\tau$ described in
\cite[Section 6.1.1]{Willbrand2005} and in Example \ref{ex.updwn}.
\end{example}

These two examples suggest
that many tests from
classical statistics  have a natural refinement by a
pre-convex rank test.
The term ``pre-convex'' refers to the following interpretation
of the axiom (PC). Consider any two vectors $u$ and $u'$
in $\R^n$, and a convex combination $u'' = \lambda u +
(1-\lambda) u'$, with $0 < \lambda < 1$.
If  $\pi, \pi', \pi''$ are the permutations of $u,u',u''$
then $\pi'' \in \mathcal{L}(\pi \cap \pi')$.
Thus the regions in $\R^n$ specified by a
 pre-convex rank test are convex cones.

\section{Convex rank tests}

A {\em fan} in $\R^n$ is a finite collection $\mathcal{F}$ of
polyhedral cones which satisfies the following properties:
(i) if $C \in \mathcal{F}$ and $C'$ is a face of $C$, then
$C' \in \mathcal{F}$,
(ii) If $C, C' \in \mathcal{F}$,
then $C \cap C'$ is a face of $C$.
Two vectors $u$ and $v $ in $\R^n$ are
{\em permutation equivalent} when $u_i < u_j$
if and only if $v_i < v_j$, and $u_i = u_j$
if and only if $v_i = v_j$  for all $i,j \in [n]$.
The permutation equivalence classes (of which there are $13$ for $n=3$) induce
a fan which we call the {\em $S_n$-fan}.
The maximal cones in the $S_n$-fan, which are the closures
of the permutation equivalence classes corresponding to total orders,
are indexed
by permutations $\delta$ in $ S_n$.
A {\em coarsening} of the $S_n$-fan is a fan $\mathcal{F}$ such that every
permutation equivalence class
of $\R^n$ is fully contained in a cone $C$ of $\mathcal{F}$;
$\mathcal{F}$ defines a partition of $S_n$ because each maximal cone of the $S_n$-fan is contained in some cone $C \in \mathcal{F}$.
We define a {\em convex rank test} to be a partition of $S_n$ defined by a coarsening of the $S_n$-fan.
We identify the fan with that test.

Two maximal cones of the $S_n$-fan
share a {\em wall} if there exists an index $k$ such that
$\delta_k = \delta'_{k+1}$, $  \delta_{k+1} = \delta'_k$ and
$\delta_i = \delta'_i$ for $i \neq k,k+1$.  That is, the corresponding permutations $\delta$ and $\delta'$ differ
by an adjacent transposition.
To such an unordered pair $\{\delta,\delta'\}$,
we associate the following conditional independence (CI) statement:
\begin{equation}
\label{CIStatement}
 \delta_k \perp \!\!\! \perp  \delta_{k+1} \,|\, \{\delta_1 , \ldots,  \delta_{k-1} \}.
 \end{equation}
This formula defines a map from the set of walls of the $S_n$-fan
 onto the set  of all CI statements
$$
\mathcal{T}_n \,\, =  \,\, \bigl\{
\, i \perp \!\!\! \perp j \,|\, K \,: \, K \subseteq [n] \backslash \{i,j\} \bigr\}. $$
The map from walls to CI statements is not injective; 
there are $(n-k-1)!(k-1)!$ walls which are labelled by the statement \eqref{CIStatement}.

\begin{figure}[thb]\label{UpDown}
\[
 \begin{xy}<15mm,0cm>:
(-.5,\sinpioverthree)  ="123"  *+!DR{123} *{\bullet};
(.5 ,\sinpioverthree)  ="132"  *+!DL{132} *{\bullet};
(1,0)                  ="312"  *+!L{312}  *{\bullet};
(.5,-\sinpioverthree)  ="321"  *+!UL{321} *{\bullet};
(-.5,-\sinpioverthree) ="231"  *+!UR{231} *{\bullet};
(-1,0)                 ="213"  *+!R{213}  *{\bullet};
   "123";"132" **@{.};
   "132";"312" **@{-};
   "312";"321" **@{.};
   "321";"231" **@{.};
   "231";"213" **@{-};
   "213";"123" **@{.};
(\sinpioverthree, 0.4)  *+!{1 \indep 3 | \emptyset} ;
(-\sinpioverthree, -0.5) *+!{1 \indep 3 | \{ 2\} } ;
 \end{xy}\]

\[
\begin{xy}
   <15mm,0mm>:
(0,0)  ="origin" ;
(\sinpioverthree, 0.5)  ="uprt"  *+!DL{1 \indep 3 | \emptyset} ;
(\sinpioverthree, -0.5)  ="dwnrt"  ;
(-\sinpioverthree, 0.5)  ="upl" ;
(-\sinpioverthree, -0.5)  ="dwnl" *+!UR{1 \indep 3 | \{ 2\} } ;
(0, 1)  ="up" ;
(0,-1)  ="dwn" ;
   "origin";"uprt" **@{--};
   "origin";"dwnrt" **@{-};
   "origin";"upl" **@{-};
   "origin";"dwnl" **@{--};
   "origin";"up" **@{-};
   "origin";"dwn" **@{-};
 \end{xy}\]
\caption{The permutohedron ${\bf P}_3$ and the $S_3$-fan projected to the plane. Each permutation is represented by
its descent vector $\delta = \delta_1 \delta_2 \delta_3$.
}
\end{figure}
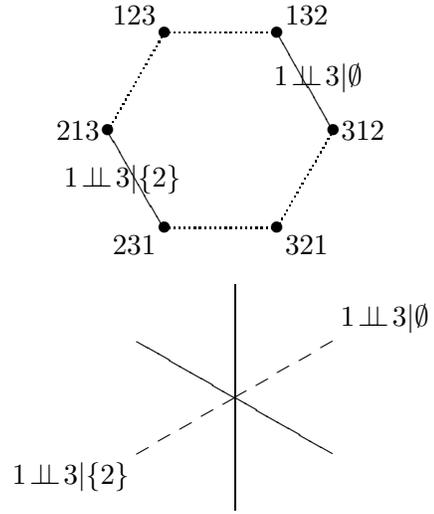

 Any convex rank test $\mathcal{F}$ is characterized
by the collection of walls $\{\delta,\delta'\}$ that are removed
when passing from the $S_n$-fan
to $\mathcal{F}$.
  So, from (\ref{CIStatement}), any convex rank test
$\mathcal{F}$ maps to a set  $\mathcal{M}_\mathcal{F}\,$
of CI statements corresponding to missing walls.
Recall from \cite{Matus2004} and \cite{Studeny2005Probabilistic}
that a subset $\mathcal{M} $ of $\mathcal{T}_n$ is a {\em
semi-graphoid} if the following axiom holds:
\begin{eqnarray*}
&  i \perp \!\!\! \perp j \, |\, K \cup {\ell}\, \in \mathcal{M}
\,\, \mbox{and} \,\,
i \perp \!\!\! \perp \ell \, |\, K\, \in \mathcal{M} \\
& \!\!\!\!\! \mbox{implies} \,\,\,
i \perp \!\!\! \perp j \,|\, K \in \mathcal{M}
\, \mbox{and} \,
i \perp \!\!\! \perp \ell \,| \, K \!\cup\! j  \in \mathcal{M}.
\end{eqnarray*}

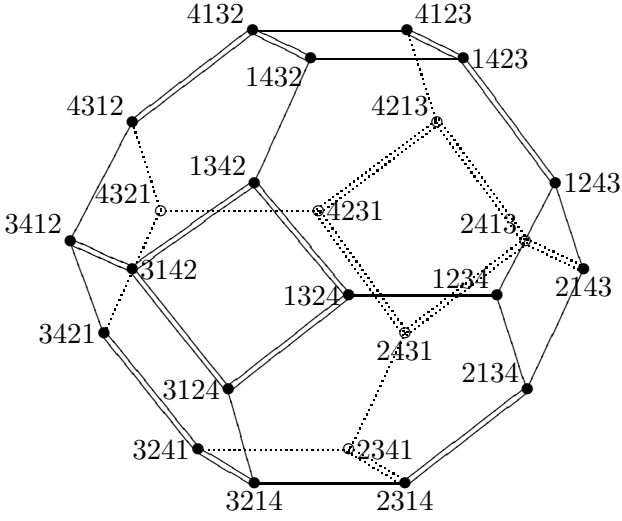
\begin{figure}[htb]\label{UpDown4}
\[
 \begin{xy}<25mm,0cm>:
%Permutohedron with n=4
%Points and labels
(1,0)  ="3214"  *+!U{3214} *{\bullet};
(1.8,0) ="2314"  *+!U{2314} *{\bullet};
(.7,.18)  ="3241"  *+!R{3241} *{\bullet};
(1.5,.18)  ="2341"  *+!L{2341} *{\circ}; %in back
(.86,.5)  ="3124"  *+!R{3124} *{\bullet};
(2.45,.5)  ="2134"  *+!DR{2134} *{\bullet};
(.2,.8)  ="3421"  *+!R{3421} *{\bullet};
(1.8 ,.8)  ="2431"  *+!U{2431} *{\circ}; %in back
(1.5,1)  ="1324"  *+!R{1324} *{\bullet};
(2.29,1)  ="1234"  *+!DR{1234} *{\bullet};
(.35,1.14)  ="3142"  *+!L{3142} *{\bullet};
(2.75,1.14)  ="2143"  *+!U{2143} *{\bullet};
(.02,1.29)  ="3412"  *+!DR{3412} *{\bullet};
(2.44,1.29)  ="2413"  *+!DR{2413} *{\circ}; %in back
(.5,1.45)  ="4321"  *+!DR{4321} *{\circ}; %in back
(1.34,1.45)  ="4231"  *+!L{4231} *{\circ}; %in back
(1,1.6)  ="1342"  *+!DR{1342} *{\bullet};
(2.6,1.6)  ="1243"  *+!L{1243} *{\bullet};
(.35,1.92)  ="4312"  *+!DR{4312} *{\bullet};
(1.97,1.92)  ="4213"  *+!DR{4213} *{\circ}; %in back
(1.3,2.26)  ="1432"  *+!UR{1432} *{\bullet};
(2.11,2.26)  ="1423"  *+!L{1423} *{\bullet};
(.99,2.41)  ="4132"  *+!DR{4132} *{\bullet};
(1.81,2.41)  ="4123"  *+!DL{4123} *{\bullet};
%EDGES %Squares
%Bottom square
"3214";"2314" **@{-}; % 2 indep 3
"3241";"2341" **@{.}; % 2 indep 3
"3241";"3214" **@{=}; % 1 indep 4 | 23
"2341";"2314" **@{:}; % 1 indep 4 | 23
%Right square
"2134";"1234" **@{-}; % 1 indep 2 |
"2143";"1243" **@{-}; % 1 indep 2 |
"1234";"1243" **@{-}; % 3 indep 4 | 12
"2134";"2143" **@{-}; % 3 indep 4 | 12
%Top square
"4132";"4123" **@{-}; % 2 indep 3 | 14
"1432";"1423" **@{-}; % 2 indep 3 | 14
"1432";"4132" **@{=}; % 1 indep 4 |
"4123";"1423" **@{=}; % 1 indep 4 |
%Left square
"4312";"3412" **@{-}; % 3 indep 4 |
"4321";"3421" **@{.}; % 3 indep 4 |
"4312";"4321" **@{.}; % 1 indep 2 | 34
"3412";"3421" **@{-}; % 1 indep 2 | 34
%Back square
"4213";"2413" **@{:}; % 2 indep 4 |
"4231";"2431" **@{:}; % 2 indep 4 |
"4213";"4231" **@{:}; % 1 indep 3 | 24
"2413";"2431" **@{:}; % 1 indep 3 | 24
%Front square
"1342";"1324" **@{=}; % 2 indep 4 | 13
"3142";"3124" **@{=}; % 2 indep 4 | 13
"1342";"3142" **@{=}; % 1 indep 3 |
"1324";"3124" **@{=}; % 1 indep 3 |
%Rest of edges
"2314";"2134" **@{=}; % 1 indep 3 | 2
"3124";"3214" **@{-}; % 1 indep 2 | 3
"3421";"3241" **@{=}; % 2 indep 4 | 3
"3412";"3142" **@{=}; % 1 indep 4 | 3
"1324";"1234" **@{-}; % 2 indep 3 | 1
"1432";"1342" **@{-}; % 3 indep 4 | 1
"4312";"4132" **@{=}; % 1 indep 3 | 4
"1423";"1243" **@{=}; % 2 indep 4 | 1
"2341";"2431" **@{.}; % 3 indep 4 | 2
"4321";"4231" **@{.}; % 2 indep 3 | 4
"2413";"2143" **@{:}; % 1 indep 4 | 2
"4123";"4213" **@{.}; % 1 indep 2 | 4
%"";"" **@{.}; %  indep  |
\end{xy}
\]
\caption{The permutohedron ${\bf P}_4$
with vertices marked by descent vectors $\delta$.
The test ``up-down analysis'' is indicated
by the double edges.}
\end{figure}

\begin{theorem} \label{fantheorem}
The map $\mathcal{F} \mapsto \mathcal{M}_\mathcal{F}$ is
a bijection between convex rank tests and  semi-graphoids.
\end{theorem}

\begin{example}[Up-down analysis for $n=3$] \label{ex.updwn}
  The test in \cite{Willbrand2005} is a convex rank test
  and is visualized in Figure 1.  Permutations are in the same class if they are connected by a solid edge; there are four classes.  In the $S_3$-fan, the two missing walls
 are labeled by conditional independence statements
 as defined in (\ref{CIStatement}).
\end{example}

\begin{example}[Up-down analysis for $n=4$]
The test $\mathcal{F}$ in \cite{Willbrand2005} is
shown in Figure 2.
The double edges correspond to the $12$ CI statements
in $\mathcal{M}_\mathcal{F}$. There are $8$ classes; e.g.,
the class $\{3412,3142,1342,1324,3124\}$ consists
of the $5$ permutations with up-down pattern  $(-,+,-)$.
\end{example}

Our proof of Theorem \ref{fantheorem} rests on translating
the semi-graphoid axiom for a set of CI statements into
geometric statements about the corresponding set of edges of the
permutohedron.

The $S_n$-fan is the normal fan \cite{Ziegler1995}
of the {\em permutohedron} ${\bf P}_n$, which is the convex hull of
the vectors $(\rho_1,\ldots,\rho_n) \in \R^n$, where $\rho$ runs
over all rank vectors of permutations in $S_n$.
The edges of ${\bf P}_n$
correspond to walls and are thus labeled with CI statements.
A collection of parallel edges of ${\bf
P}_n$ perpendicular to a hyperplane $x_i=x_j$ corresponds to the set of
CI statements $i \indep j |K$, where $K$ ranges over all subsets of
$[n] \backslash \{i,j\}$. The two-dimensional faces of ${\bf P}_n$
are squares and regular hexagons, and two edges of  ${\bf P}_n$ have
the same label in $\mathcal{T}_n$ if, but not only if, they are
opposite edges of a square.
A semi-graphoid $\mathcal{M}$ can be identified with the set $\mathbf{M}$ of
edges with labels from $\mathcal{M}$. The semi-graphoid axiom
translates into a geometric condition on the hexagonal faces of
${\bf P}_n$.

\begin{observation} \label{obs:SqHexAxioms}
A set $\mathbf{M}$ of edges of the permutohedron ${\bf P}_n$ is a
semi-graphoid if and only if  $\mathbf{M}$ satisfies the following
two axioms: \\ %\hfill \break \noindent
{\bf Square axiom:} Whenever an edge of a square is in
$\mathbf{M}$, then the opposite edge is also
in $\mathbf{M}$. \\ %\hfill \break \noindent
{\bf Hexagon axiom:} Whenever two ad\-ja\-cent edg\-es of a hexagon
are in $\mathbf{M}$, then the two opposite edges of that hexagon
are also in $\mathbf{M}$.
\end{observation}

Let $\mathbf{M}$ be the subgraph of the edge graph of ${\bf P}_n$
defined by the statements in $\mathcal{M}$. Then the classes of the rank test defined by
$\mathcal{M}$ are given by the permutations in the path-connected
components of $\mathbf{M}$. We regard a path from
$\delta$ to $\delta'$ on ${\bf P}_n$ as a word $\sigma^{(1)} \cdots
\sigma^{(l)}$ in the free associative algebra $\mathcal{A}$
generated by the adjacent transpositions of $[n]$.  For example, the
word $\sigma_{23} := (23)$ gives the path from $\delta$ to
$\delta'=\sigma_{23} \delta = \delta_1 \delta_3 \delta_2 \delta_4
\dots \delta_n$.  The following relations in $\mathcal{A}$
define a presentation of the group algebra of $S_n$:
\begin{align*}
(BS) & \; \sigma_{i, i+1} \sigma_{i+k+1, i+k+2} - \sigma_{i+k+1, i+k+2} \sigma_{i, i+1},\\
(BH) & \; \sigma_{i, i+1} \sigma_{i+1, i+2} \sigma_{i, i+1} - \sigma_{i+1, i+2} \sigma_{i, i+1} \sigma_{i+1, i+2}, %and (this word removed to fit in columns)
\\
(BN) &\;  \sigma_{i, i+1}^2 -1,
\end{align*}
where suitable $i$ and $k$ vary over $[n]$.  %where this makes sense.
The first two are the \emph{braid relations}, and the last
represents the idempotency of each transposition.

Now, we regard these relations as properties of a
set of edges of ${\bf P}_n$, by identifying a word and a
permutation $\delta$ with the set of edges that comprise the
corresponding path in ${\bf P}_n$.  For example, a set satisfying
(BS) is one such that, starting from any $\delta$, the edges of the
path $\sigma_{i, i+1} \sigma_{i+k+1, i+k+2}$ are in the set if and
only if the edges of the path $ \sigma_{i+k+1, i+k+2} \sigma_{i,
i+1}$ are in the set. Note then, that (BS) is the square axiom,
and (BH) is a weakening of the hexagon axiom of semi-graphoids.  That is, implications in
either direction hold in a semi-graphoid.  However, (BN) holds only
directionally in a semi-graphoid: if an edge lies
in the semi-graphoid, then its two vertices are in the same class;
but the empty path at some vertex $\delta$ certainly does not imply
the presence of all incident edges in the semi-graphoid.  Thus, for
a semi-graphoid, we have (BS) and (BH), but must replace (BN) with
the  directional version

\vspace{3mm}
\noindent $(BN') \;\; \qquad \qquad \sigma_{i, i+1}^2 \rightarrow 1.$
\vspace{3mm}

\noindent
Consider a path $p$ from $\delta$ to $\delta'$ in a semi-graphoid.  A result of \cite{Tits1968Problem} gives the following lemma; see also~\cite[p.~49-51]{Brown1989}.

\begin{lemma} \label{lem.allshortestpaths}
If $\mathcal{M}$ is a semi-graphoid, then if $\delta$ and $\delta'$ lie in the same class of $\mathcal{M}$, then so do all shortest paths on ${\bf P}_n$ between them.
\end{lemma}

We are now equipped to prove Theorem \ref{fantheorem}.  Note that we have demonstrated that semi-graphoids and convex rank tests can be regarded as sets of edges of ${\bf P}_n$, so we will show that their axiom systems are equivalent.  We first
show that a semi-graphoid satisfies (PC). Consider $\delta, \delta'$
in the same class $C$ of a semi-graphoid, and let $\delta'' \in
\mathcal{L}(\delta, \delta')$.  Further, let $p$ be a shortest path from
$\delta$ to $\delta''$ (so, $p \delta = \delta''$), and let $q$ be a
shortest path from $\delta''$ to $\delta'$.  We claim that $qp$ is a
shortest path from $\delta$ to $\delta'$, and thus
$\delta'' \in C$ by Lemma
\ref{lem.allshortestpaths}.  Suppose $qp$ is not a
shortest path.  Then, we can obtain a shorter path in the
semi-graphoid by some sequence of substitutions according to (BS),
(BH), and (BN').  Only (BN') decreases the length of a
path, so the sequence must involve (BN').  Therefore, there is some $i$, $j$ in $[n]$, such that their positions
relative to each other are reversed twice in $qp$. But $p$ and $q$
are shortest paths, hence  one reversal occurs in each $p$ and
$q$. Then $\delta$ and $\delta'$ agree on whether $i>j$ or
$j>i$, but the reverse holds in $\delta''$, contradicting
$\delta'' \in \mathcal{L}(\delta, \delta')$.  Thus every
semi-graphoid is a pre-convex rank test.

Now, we show that a semi-graphoid corresponds to a fan.
Consider the cone corresponding to a class $C$.  We need only
show that it meets any other cone in a shared face. Since $C$ is a cone of a coarsening of the $S_n$-fan, each nonmaximal face of $C$ lies in a hyperplane
$H =\{x_i=x_j\}$.
 Suppose a face of $C$ coincides with the hyperplane $H$
 and that $i>j$ in $C$.  A vertex $\delta$
borders $H$ if $i$ and $j$ are adjacent in $\delta$. We will
show that if $\delta,\delta' \in C$ border $H$, then their
reflections $\hat{\delta} = \delta_1 \dots ji \dots
\delta_n$ and $\hat{\delta'}= \delta'_1 \dots ji \dots \delta'_n$
both lie in some class $C'$. Consider a `great circle' path between
$\delta$ and $\delta'$ which stays closest to $H$: all vertices
in the path have $i$ and $j$ separated by at most one position, and
no two consecutive vertices have $i$ and $j$ nonadjacent.  This is a
shortest path, so it lies in $C$, by Lemma
\ref{lem.allshortestpaths}. Using the
square and hexagon axioms (Observation \ref{obs:SqHexAxioms}), we
see that the reflection of the path across $H$ is a path
in the semi-graphoid that connects $\hat{\delta}$ to
$\hat{\delta'}$ (Figure 3).  Thus a semigraphoid is a convex rank test.
\begin{figure}[htb]\label{fig:reflection}

\[
 \begin{xy}<10mm,0cm>:
%left hexagon
(0,0); %center
p+ (\sinpioverthree, 0.5) *{\bullet}; %uprt
p + (0,-1) *{\bullet} **@{.}; %dwnrt
p + (-\sinpioverthree,-.5) *{\bullet} **@{-}; %upl
p + (-\sinpioverthree,+.5) *{\bullet} *+!UR{\hat{\delta}} **@{-}; %dwnl
p + (0,1) *{\bullet} *+!DR{\delta} **@{.};
p+ (\sinpioverthree, 0.5) *{\bullet} **@{-};
p+ (\sinpioverthree, -0.5) *{\bullet} **@{-};
%connectors
(\sinpioverthree, 0.5);
p + (1,0) **@{-};
(\sinpioverthree, -0.5);
p + (1,0) **@{-};
%right hexagon
(2.73205081,0); %center
p+ (\sinpioverthree, 0.5) *+!DL{\delta'} *{\bullet}; %uprt
p + (0,-1) *{\bullet} *+!UL{\hat{\delta'}} **@{.}; %dwnrt
p + (-\sinpioverthree,-.5) *{\bullet} **@{-}; %upl
p + (-\sinpioverthree,+.5) *{\bullet} **@{-}; %dwnl
p + (0,1) *{\bullet} **@{.};
p+ (\sinpioverthree, 0.5) *{\bullet} **@{-};
p+ (\sinpioverthree, -0.5) *{\bullet} **@{-};
% Hyperplane
(-1.5,0);
p+(5.8,0) *+!UL{x_i = x_j}**@{--};
 \end{xy}\]
\caption{Reflecting a path across a hyperplane.}
\end{figure}
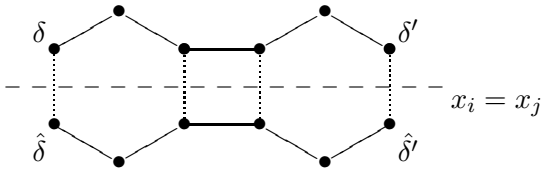

Finally, if $\mathbf{M}$ is a set of edges of ${\bf P}_n$, representing
a convex rank test, then it is easy to show that $\mathbf{M}$
satisfies the square and hexagon axioms.
This completes the proof of Theorem~\ref{fantheorem}.

\begin{remark}
For $n=3$ there are
$40$ pre-convex rank tests,
but only $22$ of them are convex rank tests.
The corresponding CI models are shown
in Figure 5.6 on page 108 in \cite{Studeny2005Probabilistic}.
\end{remark}

\section{The submodular cone}

In this section we examine a subclass of the convex rank tests.
Let $2^{[n]}$ denote the collection of all
subsets of $[n] = \{1,2,\ldots,n\}$. Any real-valued function $\, w : 2^{[n]} \rightarrow \R \, $
defines a convex polytope $Q_w$ of dimension $\leq n-1$
as follows:
\begin{eqnarray*}  Q_w & :=  \,
\bigl\{ \, x \in \R^n \,: \,
x_1 + x_2 + \cdots + x_n = w([n]) \\
& \text{\ \,and } \sum\nolimits_{i \in I} x_i \leq w(I)\,\,
\hbox{for all} \,\, \emptyset\neq I \subseteq [n]   \,\bigr\}.
\end{eqnarray*}
A function $\, w : 2^{[n]} \rightarrow \R \, $ is called
 {\em submodular} if
$\,w(I) + w(J)\, \geq\, w(I \cap J) +  w(I \cup J)\,$
for $I,J \subseteq [n]$.

\begin{proposition} \label{prop:submodularnormal}
A function $\,w: 2^{[n]} \rightarrow \R \, $
is submodular if and only if
the normal fan of the polyhedron $Q_w$
is a coarsening of the $S_n$-fan.
\end{proposition}

This  follows from
greedy maximization as in \cite{Lovasz1983Submodular}.
Note that the function $w$ is submodular if and only if
the optimal solution of
$$
\mbox{maximize $u \cdot x$ subject to $x \in Q_w$}
$$
depends only on the permutation equivalence class
of $u$.
Thus, solving this linear programming problem
constitutes a convex rank test.  Any such test is called a
{\em submodular rank test}.

A convex polytope is a {\em (Minkowski) summand}
of another polytope if the normal fan of the latter
refines the normal fan of the former. The
polytope $Q_w$ that represents a submodular rank test
is a  summand of the permutohedron
${\bf P}_n$.

\begin{theorem}
The following combinatorial objects are equivalent for any positive integer~$n$: \\
\noindent $1.$ submodular rank tests, \hfill \break
\noindent $2.$ summands of the permutohedron $\mathbf{P}_n$, \hfil \break
\noindent $3.$ structural conditional~independence~models, \hfil \break
\noindent $4.$ faces of the submodular cone ${\bf C}_n$ in $\R^{2^n}$. \hfill \break
\end{theorem}
We have 1$\iff$2 from Proposition \ref{prop:submodularnormal}, and
1$\iff$3 follows from \cite{Studeny2005Probabilistic}. Further 3$\iff$4 holds by definition.

The {\em submodular cone}
is the cone ${\bf C}_n$ of all submodular functions $w :
2^{[n]} \rightarrow \R$.
  Working modulo its lineality space
   $\,{\bf C}_n \cap (-{\bf C}_n) $, we regard
   ${\bf C}_n$ as a pointed cone of dimension $2^n-n-1$.

\begin{remark}
All $22$ convex rank tests for $n=3$ are submodular.
The submodular cone ${\bf C}_3$ is a
$4$-dimensional cone whose base is a
bipyramid. The polytopes $Q_w$, as
$w$ ranges over the faces of ${\bf C}_3$,
are all the Minkowski summands of~${\bf P}_3$.
\end{remark}

\begin{proposition}
For $n \geq 4$, there exist convex rank tests that are not submodular rank tests.
Equivalently, there are fans that coarsen the $S_n$-fan
but are not the normal fan of any polytope.
\end{proposition}

This result is stated in Section 2.2.4 of \cite{Studeny2005Probabilistic} in the following form:
 ``There exist semi-graphoids that are not structural.''

 We answered Question 4.5 posed in \cite{PRW} by finding a non-submodular convex rank test in which all the posets $P_i$ are trees:
 \begin{eqnarray*}
\mathcal{M} \! & = \,\,\,
\bigl\{
 2 \perp \!\!\! \perp 3 | \{1,4\},\,
 1 \perp \!\!\! \perp 4 |  \{2,3\}, \,  \\ &
 1 \perp \!\!\! \perp 2 | \emptyset,\,
 3 \perp \!\!\! \perp 4 |\emptyset \,\bigr\}.
 \end{eqnarray*}

\begin{remark}
For $n=4$ there are $22108$ submodular rank tests, one for each face of the
$11$-dimensional cone ${\bf C}_4$.
The base of this submodular cone is a polytope with
$f$-vector $
(1,37, 356,  $ $ 1596, 3985, 5980, 5560, 3212, 1128, 228, 24,1)$.
\end{remark}

\begin{remark}
For $n=5$ there are
$117978$ coarsest submodular rank tests,
in $1319$ symmetry classes.
We confirmed this result of \cite{Studeny2000} with {\tt POLYMAKE} \cite{Gawrilow2000}.
\end{remark}

We now define a class of submodular rank tests,
which we call {\em Minkowski sum of
simplices (MSS) tests}. Note that each subset $K$ of $[n]$
defines a submodular function $w_K$
by setting $w_K (I) = 1$ if $K \cap I $ is non-empty
and  $w_K(I) = 0$ if $K \cap I $ is empty.
The corresponding polytope
$Q_{w_K}$ is the simplex
$\Delta_K = {\rm conv} \{ e_k :  k \in K \}$.

Now consider an arbitrary subset  $\,\mathcal{K} = \{K_1,K_2,\ldots,K_r \}\,$
of $2^{[n]}$. It defines the submodular function
$\,w_{\mathcal{K}} =  w_{K_1} + w_{K_2} + \cdots + w_{K_r}$.
 The corresponding polytope is the Minkowski sum
$$ \Delta_\mathcal{K} \quad = \quad \Delta_{K_1} + \Delta_{K_2} + \cdots + \Delta_{K_r}. $$
 The associated MSS test $\tau_\mathcal{K}$ is defined as follows.
 Given $\rho \in S_n$, we compute the number of indices
 $j \in [r]$ such that $\,{\rm max}\{ \rho_k \,: \, k \in K_j \}\, = \,
\rho_i $,
 for each  $i \in [n]$.
 The signature $\tau_\mathcal{K}(\rho)$ is
the vector in $\N^n$ whose $i$th coordinate is that number. 
Few submodular rank tests are MSS tests:

\begin{remark}
For $n = 3$, among the $22$ submodular rank tests,
only $15$ are MSS tests.
For $n=4$, among the $22108$,
only $1218$ are MSS.
\end{remark}

%Section Five
\section{Graphical tests}

Graphical models are fundamental in statistics, and
 they also lead to a useful class of rank tests.
First we show how to associate a semi-graphoid to a family
$\mathcal{K}$. Let $\mathcal{F}_{w_{\mathcal{K}}}$ be the normal fan
 of $Q_{w_{\mathcal{K}}}$. We write  $\mathcal{M}_\mathcal{K}$
for the CI model derived from $\mathcal{F}_{w_{\mathcal{K}}}$ using the bijection in
Theorem \ref{fantheorem}.

\begin{proposition} \label{CIsetfam}
The semi-graphoid $\mathcal{M}_\mathcal{K}$
is the set of CI statements
$\, i \perp \!\!\! \perp j \, |\, K \,$ which
satisfy the following property:
all sets  containing $\{i,j\}$
and contained in $\{i,j\} \cup [n] \backslash K \,$
are not in $\mathcal{K}$.
\end{proposition}

Let $G$ be a graph with vertex set $[n]$.
We define $\mathcal{K}(G)$ to be the collection of all subsets $K$
of $[n]$ such that the induced subgraph of $G|_K$ is connected.
Recall that the {\em undirected graphical model}
(or {\em Markov random field}) derived from the graph $G$
is the set $\mathcal{M}^G$  of CI statements:
\begin{eqnarray*}
& \!\!\!\!\! \mathcal{M}^G \,= \,
\bigl\{\, i \perp \!\!\! \perp j \,|\,  C \,\, :\,\,
\mbox{the restriction of $G$ to} \quad  \\
& \qquad \,\,\,\,  [n] \backslash C \,\,
\mbox{ contains no path from $i$ to $j$} \bigr\}.
\end{eqnarray*}

The polytope $\Delta_G = \Delta_{\mathcal{K}(G)}$
is the {\em graph associahedron}, which is
a well-studied object in combinatorics
\cite{Carr2004,Postnikov2005}.
The next theorem is derived from Proposition~\ref{CIsetfam}.

\begin{theorem} \label{Jasonslemma}
The CI model induced by the graph associahedron
coincides with the graphical
model $\mathcal{M}^G$, i.e., $\,\mathcal{M}_{\mathcal{K}(G)}
\,= \, \mathcal{M}^G $.
\end{theorem}

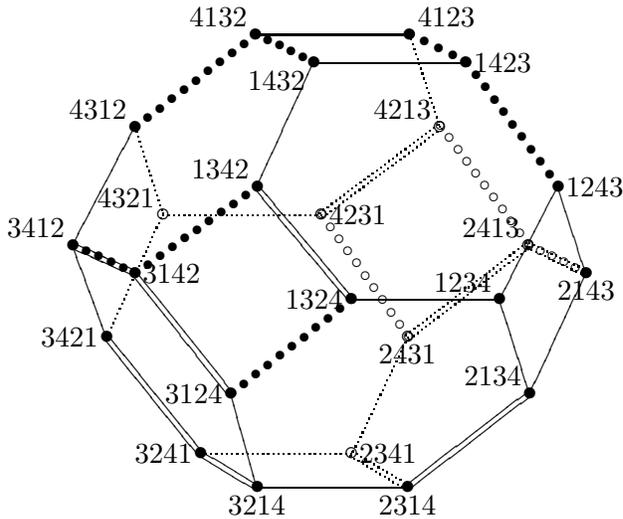
\begin{figure}[htb]\label{GM}
\[
 \begin{xy}<25mm,0cm>:
%Permutohedron with n=4
%Points and labels
(1,0)  ="3214"  *+!U{3214} *{\bullet};
(1.8,0) ="2314"  *+!U{2314} *{\bullet};
(.7,.18)  ="3241"  *+!R{3241} *{\bullet};
(1.5,.18)  ="2341"  *+!L{2341} *{\circ}; %in back
(.86,.5)  ="3124"  *+!R{3124} *{\bullet};
(2.45,.5)  ="2134"  *+!DR{2134} *{\bullet};
(.2,.8)  ="3421"  *+!R{3421} *{\bullet};
(1.8 ,.8)  ="2431"  *+!U{2431} *{\circ}; %in back
(1.5,1)  ="1324"  *+!R{1324} *{\bullet};
(2.29,1)  ="1234"  *+!DR{1234} *{\bullet};
(.35,1.14)  ="3142"  *+!L{3142} *{\bullet};
(2.75,1.14)  ="2143"  *+!U{2143} *{\bullet};
(.02,1.29)  ="3412"  *+!DR{3412} *{\bullet};
(2.44,1.29)  ="2413"  *+!DR{2413} *{\circ}; %in back
(.5,1.45)  ="4321"  *+!DR{4321} *{\circ}; %in back
(1.34,1.45)  ="4231"  *+!L{4231} *{\circ}; %in back
(1,1.6)  ="1342"  *+!DR{1342} *{\bullet};
(2.6,1.6)  ="1243"  *+!L{1243} *{\bullet};
(.35,1.92)  ="4312"  *+!DR{4312} *{\bullet};
(1.97,1.92)  ="4213"  *+!DR{4213} *{\circ}; %in back
(1.3,2.26)  ="1432"  *+!UR{1432} *{\bullet};
(2.11,2.26)  ="1423"  *+!L{1423} *{\bullet};
(.99,2.41)  ="4132"  *+!DR{4132} *{\bullet};
(1.81,2.41)  ="4123"  *+!DL{4123} *{\bullet};
%EDGES %Squares
%Bottom square
"3214";"2314" **@{-}; % 2 indep 3
"3241";"2341" **@{.}; % 2 indep 3
"3241";"3214" **@{=}; % 1 indep 4 | 23
"2341";"2314" **@{:}; % 1 indep 4 | 23
%Right square
"2134";"1234" **@{-}; % 1 indep 2 |
"2143";"1243" **@{-}; % 1 indep 2 |
"1234";"1243" **@{-}; % 3 indep 4 | 12
"2134";"2143" **@{-}; % 3 indep 4 | 12
%Top square
"4132";"4123" **@{-}; % 2 indep 3 | 14
"1432";"1423" **@{-}; % 2 indep 3 | 14
"1432";"4132" **@{*}; % 1 indep 4 |
"4123";"1423" **@{*}; % 1 indep 4 |
%Left square
"4312";"3412" **@{-}; % 3 indep 4 |
"4321";"3421" **@{.}; % 3 indep 4 |
"4312";"4321" **@{.}; % 1 indep 2 | 34
"3412";"3421" **@{-}; % 1 indep 2 | 34
%Back square
"4213";"2413" **@{o}; % 2 indep 4 |
"4231";"2431" **@{o}; % 2 indep 4 |
"4213";"4231" **@{:}; % 1 indep 3 | 24
"2413";"2431" **@{:}; % 1 indep 3 | 24
%Front square
"1342";"1324" **@{=}; % 2 indep 4 | 13
"3142";"3124" **@{=}; % 2 indep 4 | 13
"1342";"3142" **@{*}; % 1 indep 3 |
"1324";"3124" **@{*}; % 1 indep 3 |
%Rest of edges
"2314";"2134" **@{=}; % 1 indep 3 | 2
"3124";"3214" **@{-}; % 1 indep 2 | 3
"3421";"3241" **@{=}; % 2 indep 4 | 3
"3412";"3142" **@{=}; % 1 indep 4 | 3
"3412";"3142" **@{*}; % 1 indep 4 | 3
"1324";"1234" **@{-}; % 2 indep 3 | 1
"1432";"1342" **@{-}; % 3 indep 4 | 1
"4312";"4132" **@{*}; % 1 indep 3 | 4
"1423";"1243" **@{*}; % 2 indep 4 | 1
"2341";"2431" **@{.}; % 3 indep 4 | 2
"4321";"4231" **@{.}; % 2 indep 3 | 4
"2413";"2143" **@{:}; % 1 indep 4 | 2
"2413";"2143" **@{o}; % 1 indep 4 | 2
"4123";"4213" **@{.}; % 1 indep 2 | 4
%"";"" **@{.}; %  indep  |
\end{xy}
\]
\caption{The permutohedron ${\bf P}_4$. Double edges indicate the test $\tau_{\mathcal{K}(G)}$ when
$G$ is the path. Edges with large dots
indicate the test $\tau^*_{\mathcal{K}(G)}$.}
\end{figure}

There is a natural involution $*$ on the set of all CI statements
which is defined as follows:
$$ ( i \perp \!\!\! \perp j \,|\,  C)^* \quad := \quad
 i \perp \!\!\! \perp j \,|\,  [n]\backslash (C \cup \{i,j\}) . $$
If $\mathcal{M}$ is any CI model,
then the CI model $\mathcal{M}^*$ is obtained by applying the involution $*$
to all the CI statements in the model $\mathcal{M}$.  Note that this involution was called duality in
\cite{Matus1992Ascending}.  The {\em graphical
tubing rank test} $\tau^*_{\mathcal{K}(G)}$ is the test associated with $\scrM^*_{\mathcal{K}(G)}$. It can be obtained by a
construction similar to the MSS test $\tau_{\mathcal{K}}$, with
the function $w_{\mathcal{K}}$ defined differently and supermodular.
The {\em graphical model rank test} $\tau_{\mathcal{K}(G)}$ is the
 MSS test of the set family $\mathcal{K}(G)$.

We next relate $\tau_{\mathcal{K}(G)}$ and $\tau^*_{\mathcal{K}(G)}$ to a known
combinatorial characterization of
the graph associahedron $\Delta_G$.
Two subsets $A,B$ $\subset [n]$ are
\emph{compatible} for the graph $G$
if one of the following conditions holds: $A\subset B$, $B\subset A$, or $A\cap B = \emptyset$, and there is no edge between any node in
$A$ and $B$. A {\em tubing} of the graph $G$
is a subset ${\bf T}$ of $2^{[n]}$ such that
any two elements of ${\bf T}$ are compatible.
Carr and Devadoss (2005) showed that
$\Delta_G$ is a simple polytope whose
faces are in bijection with the tubings.

\begin{theorem} \label{maingraphical}
The following four
combinatorial objects are isomorphic for any graph $G$ on $[n]$: \hfill \break
\noindent $\bullet$ the graphical model rank test $\tau_{\mathcal{K}(G)}$, \hfill \break
\noindent $\bullet$ the graphical tubing rank test
$\tau^*_{\mathcal{K}(G)}$, \hfill \break
\noindent $\bullet$ the fan of the graph associahedron~$\Delta_G$, \hfill \break
\noindent $\bullet$ the simplicial complex of all tubings on $G$.

\end{theorem}

The maximal tubings of $G$ correspond to vertices of the graph associahedron
$\Delta_G$. When $G$ is the path of length $n$, then $\Delta_G$ is the
{\em associahedron}, and
when it is a cycle, $\Delta_G$ is the {\em cyclohedron}. The number of classes in the
tubing test $\tau^*_{\mathcal{K}(G)}$ is the $G$-Catalan number of
\cite{Postnikov2005}.  This number is
 $\frac{1}{n+1} {2n \choose n}$ for the {\em associahedron test}
 and ${2n-2 \choose n-1}$ for the {\em cyclohedron test}.

\section{Enumerating linear extensions}

In this paper we introduced a hierarchy of rank tests, ranging from
pre-convex to graphical. Rank tests are applied to data
vectors $u \in \R^n$, or permutations $\pi \in S_n$, and locate
their cones. In order to determine the significance of a data
vector, one needs to compute the quantity $\,|\, \tau^{-1} \bigl(
\tau(\pi) \bigr)\,|$, and possibly the probabilities of other
maximal cones. These cones are indexed by posets $P_1,P_2,\ldots,
P_k$ on $[n]$, and the probability computations are equivalent to
finding the cardinality of some of the sets $\mathcal{L}(P_i)$.

We now present our methods for computing linear extensions. If the
rank test is a tubing test then this computation is done as follows.
From the given permutation, we identify its signature (image under $\tau$), which we may assume is its $G$-tree ${\bf T}$ \cite{Postnikov2005}. Suppose the root of the tree ${\bf T}$ has $k$
children, each of which is a root of a subtree ${\bf T}^i$ for
$i=1,\ldots,k$. Writing $|{\bf T}^i |$ for the number of nodes in
${\bf T}^i$, we have
$$ |\, \tau^{-1}({\bf T} ) \, |
= \binom{\sum_{i=1}^k |{\bf T}^i |}{ |{\bf T}^1|, \ldots, |{\bf
T}^k|} \left( \prod_{i=1}^k |\tau^{-1}( {\bf T}^{i})| \right).$$
This recursive formula can be translated into an efficient iterative
algorithm. In \cite{Willbrand2005} the analogous problem is raised
for the test in Example 3. A determinantal formula for
(\ref{Pvalue}) appears in \cite[page 69]{Stanley1997}.

For an arbitrary convex rank test we proceed as follows.
The test is specified (implicitly or explicitly)
 by a collection of posets  $P_1,\ldots,P_k$ on $[n]$.
From the given permutation, we first identify the unique poset $P_i$
of which that permutation is a linear extension. We next construct
the {\em distributive lattice} $L(P_i)$ of all order ideals of
$P_i$. Recall that an {\em order ideal} is  a subset $O$ of $[n]$
such that if $l \in O$ and $(k,l) \in P_i$ then $k \in O$. The set
of all order ideals is a lattice with meet and join operations given
by set intersection $O \cap O'$ and set union $O \cup O'$. Knowledge
of this distributive lattice $L(P_i)$ solves our problem because the
linear extensions of $P_i$ are precisely the maximal chains of
$L(P_i)$. 
Computing the number of linear
extensions is \#P-complete \cite{Brightwell1991}. Therefore we
developed efficient heuristics to build $L(P_i)$.

The key algorithmic task is the following: given a poset $P_i$
on $[n]$, compute an efficient representation of the distributive
lattice $L(P_i)$.
Our program for performing rank tests works as follows.
 The input is a permutation $\pi$ and a rank test $\tau$.
  The test $\tau$ can be specified either
 \begin{itemize}
 \item by a list of posets $P_1,\ldots,P_k$ (pre-convex),
 \item or by a semigraphoid $\mathcal{M}$ (convex rank test),
 \item or by a submodular function $w : 2^{[n]} \rightarrow \R$,
 \item or by a collection $\mathcal{K}$ of subsets of $[n]$ (MSS),
\item or by a graph $G$ on $[n]$ \  (graphical test).
 \end{itemize}
 The output of our program has two parts.
 First, it gives the number $|\mathcal{L}(P_i)|$ of linear extensions,
  where the poset $P_i$ represents the equivalence
 class of $S_n$ specified by the data $\pi$.
It also gives a representation of
 the distributive lattice $L(P_i)$, in a format
 that can be read by the {\bf maple} package
 {\bf posets}  \cite{Stembridge2004}.
  Our software for the above rank tests is available at
$\, {\tt www.bio.math.berkeley.edu/ranktests/} $.

\medskip

\section*{Acknowledgments}

This paper originated in
discussions with Olivier Pourqui\'{e} and Mary-Lee Dequ\'{e}ant in
the DARPA Fundamental Laws of Biology Program, which supported Jason Morton, Lior Pachter, and Bernd Sturmfels.  Anne Shiu was supported by a Lucent Technologies Bell Labs Graduate Research Fellowship.  Oliver Wienand was supported by the Wipprecht foundation.


\begin{thebibliography}{}

\bibitem[\protect\citename{Brightwell}1991]{Brightwell1991}
G~Brightwell and P~Winkler.
\newblock Counting linear extensions.
\newblock {\em Order}, 8(3):225-242, 1991.

\bibitem[\protect\citename{Brown}1989]{Brown1989}
%\bibitem{Brown1989}
K~Brown.
\newblock {\em Buildings}.
\newblock Springer, New York, 1989.

\bibitem[\protect\citename{Carr}2004]{Carr2004}
%\bibitem{Carr2004}
M~Carr and S~Devadoss.
\newblock Coxeter complexes and graph associahedra.
\newblock 2004.
\newblock Available from {\tt http://arxiv.org/abs/math.QA/0407229}.

\bibitem[\protect\citename{Gawrilow}2000]{Gawrilow2000}
%\bibitem{Gawrilow2000}
E~Gawrilow and M~Joswig.
\newblock Polymake: a framework for analyzing convex polytopes, in
\newblock {\em Polytopes -- Combinatorics and Computation}, eds. G Kalai and G M Ziegler, Birkh\"auser, 2000, 43-74.


\bibitem[\protect\citename{Lov\'asz}1983]{Lovasz1983Submodular}
L~Lov\'asz.
\newblock Submodular functions and convexity, in
\newblock {\em Math Programming: The State of the Art}, eds. A~Bachem, M~Groetschel, and B~Korte, Springer, 1983, 235-257.


\bibitem[\protect\citename{Mat\'{u}\v{s}}1992]{Matus1992Ascending}
%\bibitem{Matus1992}
F~Mat\'{u}\v{s}.
\newblock Ascending and descending conditional independence relations, in
\newblock {\em Proceedings of the Eleventh Prague Conference on Inform. Theory, Stat. Dec. Functions and Random Proc.}, Academia, B, 1992, 189-200.

\bibitem[\protect\citename{Mat\'{u}\v{s}}2004]{Matus2004}
F~Mat\'{u}\v{s}.
\newblock Towards classification of semigraphoids.
\newblock {\em Discrete Mathematics}, 277, 115-145, 2004.

\bibitem[\protect\citename{Pitman}1937]{Pitman1937SignificanceI}
%\bibitem{Pitman1937SignificanceI}
EJG~Pitman.
\newblock Significance tests which may be applied to samples from any
  populations.
\newblock {\em Supplement to the Journal of the Royal Statistical Society},
  4(1):119--130, 1937.

\bibitem[\protect\citename{Postnikov}2005]{Postnikov2005}
%\bibitem{Postnikov2005}
A~Postnikov.
\newblock Permutohedra, associahedra, and beyond.
\newblock 2005.
\newblock Available from {\tt http://arxiv.org/abs/math/0507163}.

\bibitem[\protect\citename{Postnikov}2006]{PRW}
%\bibitem{PRW}
A~Postnikov, V~Reiner, L~Williams.
\newblock Faces of Simple Generalized Permutohedra.
\newblock Preprint, 2006.
%Not on the arXiv.

% we commented out the para that refers to this
%\bibitem[\protect\citename{Reading}2006]{Reading2006}
%\bibitem{Reading2006}
%N~Reading.
%\newblock Cambrian lattices.
%\newblock {\em Advances in Mathematics}, to appear.

\bibitem[\protect\citename{Stanley}1997]{Stanley1997}
%\bibitem{Stanley1997}
RP~Stanley.
\newblock {\em Enumerative Combinatorics} Volume I, Cambridge University
Press, Cambridge, 1997.


\bibitem[\protect\citename{Stembridge}2004]{Stembridge2004}
%\bibitem{Stembridge2004}
J~Stembridge.
\newblock Maple packages for symmetric functions, posets, root systems, and finite Coxeter groups.
\newblock Available from {\tt www.math.lsa.umich.edu/$\sim$jrs/maple.html}.


\bibitem[\protect\citename{Studen\'y}2005]{Studeny2005Probabilistic}
%\bibitem{Studeny2005Probabilistic}
M~Studen\'{y}.
\newblock {\em Probablistic conditional independence structures}.
\newblock Springer Series in Information Science and Statistics, Springer-Verlag, London, 2005.


\bibitem[\protect\citename{Studen\'y}2000]{Studeny2000}
%\bibitem{Studeny2000}
M~Studen\'{y}, RR~Bouckaert, and T~Kocka.
\newblock Extreme supermodular set functions over five variables.
\newblock {\em Institute of Information Theory and Automation}, Research report n. 1977, Prague, 2000.

\bibitem[\protect\citename{Tits}1968]{Tits1968Problem}
%\bibitem{Tits1968Problem}
J~Tits.
\newblock {\em Le probl\`eme des mots dans les groupes de Coxeter}.
\newblock Symposia Math., 1:175-185, 1968.


\bibitem[\protect\citename{Willbrand}2005]{Willbrand2005}
%\bibitem{Willbrand2005}
K~Willbrand, F~Radvanyi, JP~Nadal, JP~Thiery, and T~Fink.
\newblock Identifying genes from up-down properties of microarray expression
  series.
\newblock {\em Bioinformatics}, 21(20):3859--3864, 2005.

\bibitem[\protect\citename{Ziegler}1995]{Ziegler1995}
%\bibitem{Ziegler1995}
G~Ziegler.
\newblock {\em Lectures on polytopes.}
\newblock Vol. 152 of Graduate Texts in Mathematics.
\newblock Springer-Verlag, 1995.

\end{thebibliography}
\end{document}